\documentclass[11pt]{article}

\usepackage{latexsym}
\usepackage{amssymb}

\newtheorem{theorem}{Theorem}
\newtheorem{proposition}{Proposition}
\newtheorem{corollary}{Corollary}
\newtheorem{definition}{Definition}

\begin{document}
{

\begin{center}
{\large\bf 
The truncated multidimensional moment problem: canonical solutions.}
\end{center}
\begin{center}

{\bf S.M. Zagorodnyuk}
\end{center}

\noindent
\textbf{Abstract.} 
For the truncated multidimensional moment problem we introduce a notion of a canonical solution.
Namely, canonical solutions are those solutions which are generated by commuting self-adjoint extensions inside the associated Hilbert space.
It is constructed a 1-1 correspondence between canonical solutions and flat extensions of the given moments
(both sets may be empty).
In the case of the two-dimensional moment problem (with triangular truncations) a search for canonical solutions leads to
an algebraic system of equations. 
A notion of the index $i_s$ of nonself-adjointness for a set of prescribed moments is introduced.
The case $i_s=0$ corresponds to flatness. 
In the case $i_s=1$ we get explicit necessary and sufficient conditions for the existence
of canonical solutions. 
These conditions are valid for arbitrary sizes of truncations.
In the case $i_s=2$ we get either explicit conditions for the existence
of canonical solutions or a single quadratic equation with several unknowns.
Numerical examples are provided.



\section{Introduction.}

We shall need some notations in what follows. As usual, we denote by $\mathbb{R}, \mathbb{C}, \mathbb{N}, \mathbb{Z}, \mathbb{Z}_+$
the sets of real numbers, complex numbers, positive integers, integers and non-negative integers,
respectively. Let $n\in\mathbb{N}$. 
By $\mathbb{Z}^n_+$ we mean $\mathbb{Z}_+\times \ldots \times\mathbb{Z}_+$, and $\mathbb{R}^n = \mathbb{R}\times \ldots \times\mathbb{R}$,
where the Cartesian products are taken with $n$ copies.
Let $\mathbf{k} = (k_1,\ldots,k_n)\in\mathbb{Z}^n_+$, $\mathbf{t} = (t_1,\ldots,t_n)\in\mathbb{R}^n$. Then
$\mathbf{t}^{\mathbf{k}}$ means the monomial $t_1^{k_1}\ldots t_n^{k_n}$, and $|\mathbf{k}| = k_1 + \ldots + k_n$.
By $\mathfrak{B}(\mathbb{R}^n)$ we denote the set of all Borel subsets of $\mathbb{R}^n$.

Let $\mathcal{K}$ be an arbitrary finite subset of $\mathbb{Z}^n_+$. Let $\mathcal{S} = (s_{\mathbf{k}})_{\mathbf{k}\in\mathcal{K}}$
be an arbitrary set of real numbers.
\textit{The truncated multidimensional moment problem} consists of finding a (non-negative) measure $\mu$ on $\mathfrak{B}(\mathbb{R}^n)$
such that
\begin{equation}
\label{f1_1}
\int \mathbf{t}^{\mathbf{k}} d\mu(\mathbf{t}) = s_{\mathbf{k}},\qquad  \forall \mathbf{k}\in\mathcal{K}.
\end{equation}

The multidimensional moment problems were described in the books of Berezansky~\cite{cit_1000_Berezansky_1965__Book},
Berg, Christensen and Ressel~\cite{cit_980_Berg_Christiansen_Ressel__Book}, 
Curto and Fialkow~\cite{cit_985_Curto_Fialkow__Book1}, \cite{cit_985_Curto_Fialkow__Book2},
Schm\"udgen~\cite{cit_1000_Schmudgen_Book_2017}, 
Shohat and Tamarkin~\cite{cit_1000_Shohat_Tamarkin_1943__Book}.
Besides the well investigated case $n=1$,
much success can be achieved in the case of (the full) $K$-moment problems, i.e. when the support of $\mu$ is restricted to a suitable 
(algebraic or semi-algebraic) prescribed set $K$. 
Here simple conditions for the solvability were obtained by
Schm\"udgen, see a detailed exposition of the corresponding results in~\cite{cit_1000_Schmudgen_Book_2017}.
In some special cases an analytical parametrization of all solutions is available, see~\cite{cit_14700_Zagorodnyuk_2013_MFAT} for
the case of a strip.
When there are no additional assumptions on the support of $\mu$, multidimentional moment problems with $n\geq 2$
appeared to be very hard for investigation.
In the cases of low numbers of prescribed moments and some restrictions on the measure support
some satisfactory results for truncated multidimentional moment problems were obtained earlier, see, e.g.,
\cite{cit_987_Fialkow_2011}, \cite{cit_14000_Yoo}, \cite{cit_15000_Zagorodnyuk_2018_AOT} and references therein.
In the case of arbitrary sizes of truncations only flatness and few other special conditions can be mentioned, 
see~\cite{cit_985_Curto_Fialkow__Book1}, \cite{cit_985_Curto_Fialkow__Book2}.
Some complicated conditions for the solvability of multidimentional moment problems were given 
in~\cite{cit_1000_Shohat_Tamarkin_1943__Book}, \cite{cit_1100_Stochel_Szafraniec__JFA_1998},
\cite{cit_993_P_V__1999}, \cite{cit_14500_Zagorodnyuk_2010_AFA}, 
\cite{cit_982_C_St_Sz_2011}, \cite{cit_14220_Vasilescu}.
Instead of restricting the measure support in~(\ref{f1_1}) one can extend the support from $\mathbb{R}^n$ to $\mathbb{C}^n$.
In this case simple conditions for the solvability are available, see~\cite{cit_3540_Zagorodnyuk_2022_Axioms}.

In this paper we shall focus on the truncated multidimensional moment problem.
We shall use the operator approach to this moment problem, as it was described 
in~\cite{cit_20000_Zagorodnyuk_2019_Concr_Oper}.
Let us recall some related definitions and notations.

Consider the following operator $W_j$ on $\mathbf{Z}^n_+$, which increases the $j$-th coordinate:
\begin{equation}
\label{f2_1}
W_j (k_1, \ldots, k_{j-1}, k_j, k_{j+1},\ldots, k_n) = (k_1, \ldots, k_{j-1}, k_j + 1, k_{j+1},\ldots, k_n),
\end{equation}
for $j=1,\ldots,n$.
A finite subset $K\subset \mathbb{Z}^n_+$ is said to be \textit{admissible}, if the following conditions hold:

\begin{itemize}
\item[1)] $\mathbf{0} = (0,\ldots,0)\in K$;

\item[2)] $\forall \mathbf{k}\in K\backslash\{ \mathbf{0} \}$, 
\begin{equation}
\label{f2_5}
\mathbf{k} = W_{a_{|\mathbf{k}|}} W_{a_{|\mathbf{k}| - 1}} \ldots W_{a_1} \mathbf{0},
\end{equation}
for some $a_j\in\{ 1,\ldots,n \}$, and
\begin{equation}
\label{f2_7}
\widetilde{\mathbf{k}}_r := W_{a_r} \ldots W_{a_1} \mathbf{0} \in K,\qquad \forall r=1,2,\ldots, |\mathbf{k}|.
\end{equation}
\end{itemize}

Such truncations appeared earlier in a paper by Laurent and Mourrain~\cite{cit_988_Laurent_Mourrain}.
There are two important special cases. The first one is called \textit{the case of triangular truncations}: 
\begin{equation}
\label{pf2_7_5}
K = K_r = \{ \mathbf{k}\in\mathbb{Z}^n_+: |\mathbf{k}|\leq r \},\ r\in\mathbb{Z}_+. 
\end{equation}
The second one is called \textit{the case of rectangular truncations}:
$$ K = K_{d_1,d_2,\ldots,d_n} = \{ \mathbf{k}=(k_1,\ldots,k_n)\in\mathbb{Z}^n_+: k_1\leq d_1,\ldots,k_n\leq d_n \}, $$ 
$$ d_1,\ldots, d_n\in\mathbb{Z}_+. $$ 
For practical purposes, it is convenient to introduce some
indexation in the set $K$ by a unique index $j$:
\begin{equation}
\label{f2_9}
K = \left\{
\mathbf{k}_j(\in\mathbb{Z}^n_+),\quad j=0,1,\ldots,\rho 
\right\}.
\end{equation}
Set
$\vec e_r := (\delta_{r,m})_{m=1}^n \in \mathbb{Z}^n_+$, $r=1,\ldots,n$, and
\begin{equation}
\label{f2_31}
\Omega_l = \{ j\in \{ 0,\ldots,\rho \}:\ \mathbf{k}_j + \vec e_l \in K \},\qquad l=1,\ldots,n.
\end{equation}
We also denote
\begin{equation}
\label{f2_23}
\Gamma = \left( s_{\mathbf{k}_j + \mathbf{k}_m} \right)_{m,j=0}^\rho,
\end{equation}
\begin{equation}
\label{f2_39}
\Gamma_l = \left( s_{\mathbf{k}_j + \mathbf{k}_m} \right)_{m,j\in\Omega_l},\quad
\widehat \Gamma_l = \left( s_{\mathbf{k}_j + \vec e_l + \mathbf{k}_m + \vec e_l} \right)_{m,j\in\Omega_l},\qquad l=1,2,\ldots,n,
\end{equation}
where the indices from $\Omega_l$ are taken in the increasing order.
Then we have the following necessary conditions of the solvability:
\begin{equation}
\label{f2_25}
\Gamma \geq 0,
\end{equation}
\begin{equation}
\label{f2_45}
\mathop{\rm Ker}\nolimits \Gamma_l \subseteq \mathop{\rm Ker}\nolimits \widehat \Gamma_l,\qquad l=1,2,\ldots,n.
\end{equation}
We do not know, if conditions~(\ref{f2_25}),(\ref{f2_45}) are sufficient for the solvability of
the moment problem~(\ref{f1_1}). This is an open problem.

Now suppose that for an admissible finite set $K\subset\mathbb{Z}^n_+$ the moment problem~(\ref{f1_1}), with
$\mathcal{K} = K+K$ and some $\mathcal{S} = (s_{\mathbf{k}})_{\mathbf{k}\in\mathcal{K}}$, is given.
Choose and fix some indexation~(\ref{f2_9}). Assume that conditions~(\ref{f2_25}),(\ref{f2_45})
hold.
We may construct \textit{the associated Hilbert space} $H$ in the following way. A set $\mathfrak{L}$ of all polynomials of the following form:
\begin{equation}
\label{f2_20}
p(\mathbf{t}) = \sum_{j=0}^\rho \alpha_j \mathbf{t}^{\mathbf{k}_j},\qquad \alpha_j\in\mathbb{C},
\end{equation}
is a linear vector space. Consider the following functional:
$$ <p, q> = \sum_{j,m=0}^\rho \alpha_j \overline{\beta_m} s_{\mathbf{k}_j + \mathbf{k}_m}, $$
where $p$ is from~(\ref{f2_20}), and $q$ has the same form as $p$, but with $\beta_j (\in\mathbb{C})$ instead of $\alpha_j$.
Introducing the classes of the equivalence $[p] = [p]_{\mathfrak{L}}$ in $\mathfrak{L}$ we obtain a finite-dimensional Hilbert space $H$.
For $l=1,\ldots,n$ we consider the following operators:
\begin{equation}
\label{f3_5}
M_l \sum_{j\in\Omega_l} \alpha_j [ \mathbf{t}^{\mathbf{k}_j} ]_{\mathfrak{L}} =
\sum_{j\in\Omega_l} \alpha_j [ \mathbf{t}^{\mathbf{k}_j + \vec e_l} ]_{\mathfrak{L}},\quad \alpha_j\in\mathbb{C}, 
\end{equation}
with $D(M_l) = \mathop{\rm Lin}\nolimits \{ [ \mathbf{t}^{\mathbf{k}_j} ]_{\mathfrak{L}} \}_{j\in\Omega_l}$.
In particular, we have
\begin{equation}
\label{f3_7}
M_l [ \mathbf{t}^{\mathbf{k}_j} ]_{\mathfrak{L}} =
[ \mathbf{t}^{\mathbf{k}_j + \vec e_l} ]_{\mathfrak{L}},\qquad j\in\Omega_l;\quad l=1,\ldots,n. 
\end{equation}
Operators $M_l$ are symmetric and they may be defined on proper subspaces of $H$.
Denote
\begin{equation}
\label{f3_30}
\Omega_0 = \{ j\in \{0,\ldots,\rho \}:\ \mathbf{k}_j + \vec e_1, \mathbf{k}_j + \vec e_2, \ldots, \mathbf{k}_j + \vec e_n \in K \},
\end{equation}
and
\begin{equation}
\label{f3_35}
H_0 = \mathop{\rm Lin}\nolimits \{ [ \mathbf{t}^{\mathbf{k}_j} ]_{\mathfrak{L}} \}_{j\in\Omega_0}.
\end{equation}
The set of moments $\mathcal{S}$ is said to be \textit{dimensionally stable}, if
$\dim H = \dim H_0$.
The set of moments $\mathcal{S}$ is said to be \textit{completely self-adjoint}, if
the operators $M_l$ are self-adjoint and pairwise commute.
Our first aim here will be to show that the dimensional stability implies the complete self-adjointness for general types of
truncations (Theorem~1).
In the case of triangular truncations it was proved by Vasilescu that flatness is equivalent to dimensional stability, see~\cite{cit_14200_Vasilescu}.
The dimensional stability can be checked numerically (see Theorem~5 in~\cite{cit_20000_Zagorodnyuk_2019_Concr_Oper})
by solving linear algebraic systems of equations and then checking some equalities.

\begin{definition}
\label{d1_1}
A solution $\mu$ of the moment problem~(\ref{f1_1}) is said to be \textbf{canonical} if it is generated by 
commuting self-adjoint operators $\widetilde M_j\supseteq M_j$ ($j=1,\ldots,n$) in the associated Hilbert space $H$, 
as it was described in the proof
of Proposition~1 in~\cite{cit_20000_Zagorodnyuk_2019_Concr_Oper}.
\end{definition}

Our second aim here is to establish a bijection between all canonical solutions of the moment problem~(\ref{f1_1}) with
general truncations
and 
all dimensionally stable close extensions
of $\mathcal{S}$ (Theorem~2). In particular, in the case of triangular truncations we get a 1-1 correspondence between
all canonical solutions and all flat extensions of $\mathcal{S}$.

In Section~3 we search for canonical solutions of the two-dimensional truncated moment problem in the case of triangular truncations. 
This leads to a system of algebraic equations.
Precise answers on the existence of canonical solutions can be given for some special cases close to flatness.
The number
\begin{equation}
\label{pf1_5}
i_s := \mathrm{dim}\ (H\ominus H_0),
\end{equation}
is said to be \textbf{the index of nonself-adjointness of $\mathcal{S}$}.
It is clear that the case $i_s=0$ means the dimension stability or flatness of $\mathcal{S}$. 
In the case $i_s=1$ we obtain explicit necessary and sufficient conditions for the existence
of canonical solutions. 
These conditions are valid for arbitrary sizes of truncations.
In the case $i_s=2$ we get either explicit conditions for the existence
of canonical solutions or a single quadratic equation with several unknowns.
Numerical examples for these cases are provided.

\noindent
{\bf Notations. }
Besides the given above notations we shall use the following conventions.
By $\mathbb{Z}_{k,l}$ we mean all integers $r$, which satisfy the following inequality:
$k\leq r\leq l$.
By $\mathop{\rm Tr}\nolimits M$ we denote the trace of a square complex matrix $M$.
If H is a Hilbert space then $(\cdot,\cdot)_H$ and $\| \cdot \|_H$ mean
the scalar product and the norm in $H$, respectively.
Indices may be omitted in obvious cases.
For a linear operator $A$ in $H$, we denote by $D(A)$
its  domain, by $R(A)$ its range, and $A^*$ means the adjoint operator
if it exists. If $A$ is invertible then $A^{-1}$ means its
inverse. $\overline{A}$ means the closure of the operator, if the
operator is closable. If $A$ is bounded then $\| A \|$ denotes its
norm.
For a set $M\subseteq H$
we denote by $\overline{M}$ the closure of $M$ in the norm of $H$.
By $\mathop{\rm Lin}\nolimits M$ we mean
the set of all linear combinations of elements from $M$,
and $\mathop{\rm span}\nolimits M:= \overline{ \mathop{\rm Lin}\nolimits M }$.
By $E_H$ we denote the identity operator in $H$, i.e. $E_H x = x$,
$x\in H$. In obvious cases we may omit the index $H$. If $H_1$ is a subspace of $H$, then $P_{H_1} =
P_{H_1}^{H}$ is an operator of the orthogonal projection on $H_1$
in $H$.

\section{Canonical solutions and flat extensions.}

Suppose that for an admissible finite set $K\subset\mathbb{Z}^n_+$ the moment problem~(\ref{f1_1}), with
$\mathcal{K} = K+K$ and some $\mathcal{S} = (s_{\mathbf{k}})_{\mathbf{k}\in\mathcal{K}}$, is given.
Choose and fix some indexation~(\ref{f2_9}). Assume that conditions~(\ref{f2_25}),(\ref{f2_45})
are satisfied.
Construct the associated Hilbert space $H$ and symmetric operators $M_k$, as it was described in the Introduction.
In the case of arbitrary admissible truncations, a notion of flatness was introduced by
Laurent and Mourrain~\cite{cit_988_Laurent_Mourrain}.
It requires a use of a closure of a smaller set of monomials.
However, two different sets of monomials may have the same closure, e.g., consider truncations $K_3$ and $K_3\backslash\{ (1,1) \}$ for $n=2$.
These causes some difficulties.
Here we shall mainly use the notion of dimensional stability when dealing with general truncations.

\begin{theorem}
\label{t2_1}
Let $n\in\mathbb{N}$, $K\subset\mathbb{Z}^n_+$ be an arbitrary admissible finite set, and $\mathcal{K} = K+K$.
Let the moment problem~(\ref{f1_1}) with some prescribed moments $\mathcal{S} = (s_{\mathbf{k}})_{\mathbf{k}\in\mathcal{K}}$ be given,
and conditions~(\ref{f2_25}),(\ref{f2_45}) be satisfied (with some indexation of $K$).
If the set of moments $\mathcal{S}$ is dimensionally stable then $\mathcal{S}$ is completely self-adjoint. 
In the case of triangular truncations the converse implication is also true: if $\mathcal{S}$ is completely self-adjoint
then $\mathcal{S}$ is dimensionally stable.

\end{theorem}
\textbf{Proof.}
Let us prove the first statement of the theorem.
Here we shall use some ideas of Vasilescu from~\cite{cit_14200_Vasilescu}, but our proof is not the same.
Suppose that $\mathcal{S}$ is dimensionally stable. Since $H$ has a finite dimension, then $H=H_0$. 
Thus, all multiplication operators $M_k$ are defined on the whole space $H$. 
It remains to prove that they commute.
Choose arbitrary $a,b\in 1,...,n$, $a\not= b$. For an arbitrary $[p(\mathbf{t})]\in H_0 (=H)$,
\begin{equation}
\label{pf2_20}
p(\mathbf{t}) = \sum_{j\in\Omega_0} c_j \mathbf{t}^{\mathbf{k}_j},\qquad c_j\in\mathbb{C},
\end{equation}
we may write:
\begin{equation}
\label{pf2_22}
M_a p(\mathbf{t}) = \left[
t_a p(\mathbf{t})
\right] = [g(\mathbf{t})],
\end{equation}
where
$[g(\mathbf{t})]\in H_0$,
\begin{equation}
\label{pf2_24}
g(\mathbf{t}) = \sum_{j\in\Omega_0} \xi_j \mathbf{t}^{\mathbf{k}_j},\qquad \xi_j\in\mathbb{C}.
\end{equation}
The last equality in~(\ref{pf2_22}) follows from $H=H_0$ and the definition of $H_0$.
On the other hand, we have
\begin{equation}
\label{pf2_26}
M_b p(\mathbf{t}) = \left[
t_b p(\mathbf{t})
\right] = [h(\mathbf{t})],
\end{equation}
where
$[h(\mathbf{t})]\in H_0$,
\begin{equation}
\label{pf2_28}
h(\mathbf{t}) = \sum_{j\in\Omega_0} \gamma_j \mathbf{t}^{\mathbf{k}_j},\qquad \gamma_j\in\mathbb{C}.
\end{equation}
Then
\begin{equation}
\label{pf2_30}
M_b M_a [ p(\mathbf{t}) ] = [t_b g(\mathbf{t})],
\end{equation}
and
\begin{equation}
\label{pf2_32}
M_a M_b [ p(\mathbf{t}) ] = [t_a h(\mathbf{t})].
\end{equation}
For an arbitrary $[u(\mathbf{t})]\in H_0 (=H)$,
\begin{equation}
\label{pf2_34}
u(\mathbf{t}) = \sum_{m\in\Omega_0} d_m \mathbf{t}^{\mathbf{k}_m},\qquad d_m\in\mathbb{C},
\end{equation}
we may write:
$$ ( M_b M_a [p], [u] )_H =  ( [t_b g(\mathbf{t})], [u(\mathbf{t})] )_H = ( [g(\mathbf{t})], [t_b u(\mathbf{t})] )_H = $$
$$ = ( [t_a p(\mathbf{t})], [t_b u(\mathbf{t})] )_H = \sum_{j,m\in\Omega_0} c_j \overline{d_m} 
( [ \mathbf{t}^{ \mathbf{k}_j + \vec e_a } ], [ \mathbf{t}^{ \mathbf{k}_m + \vec e_b } ] )_H = $$
$$ = \sum_{j,m\in\Omega_0} c_j \overline{d_m} s_{\mathbf{k}_j + \mathbf{k}_m + \vec e_a + \vec e_b}; $$
$$ ( M_a M_b [p], [u] )_H =  ( [t_a h(\mathbf{t})], [u(\mathbf{t})] )_H = ( [h(\mathbf{t})], [t_a u(\mathbf{t})] )_H = $$
$$ = ( [t_b p(\mathbf{t})], [t_a u(\mathbf{t})] )_H = \sum_{j,m\in\Omega_0} c_j \overline{d_m} 
( [ \mathbf{t}^{ \mathbf{k}_j + \vec e_b } ], [ \mathbf{t}^{ \mathbf{k}_m + \vec e_a } ] )_H = $$
$$ = \sum_{j,m\in\Omega_0} c_j \overline{d_m} s_{\mathbf{k}_j + \mathbf{k}_m + \vec e_a + \vec e_b}. $$
Therefore $M_b M_a [p] = M_a M_b [p]$. Thus, operators $M_k$ are commuting self-adjoint operators.

Now assume that we have the case of triangular truncations and $\mathcal{S}$ is completely self-adjoint. 
Then the associated multiplication operators $M_k$ are self-adjoint. Thus they are defined on the whole $H$:
$D(M_k)=H$, $k=1,...,n$. Notice that for triangular truncations we have $D(M_k)=H_0$. Then $H=H_0$, and
the dimensional stability follows.
$\Box$

It is not known for the case of rectangular truncations if complete self-adjointness implies dimensional stability. It is 
an interesting open problem. Of course, the general case of arbitrary admissible truncations is also interesting for further investigations.

We shall now study canonical solutions of the moment problem~(\ref{f1_1}). 
As it was defined in the Introduction,
a solution $\mu$ of the moment problem~(\ref{f1_1}) is called canonical if it is generated by 
some commuting self-adjoint operators $\widetilde M_j\supseteq M_j$ ($j=1,\ldots,n$) in $H$.

\begin{definition}
\label{d2_1}
Let $K\subset\mathbb{Z}^n_+$ be an arbitrary admissible finite set; $n\in\mathbb{N}$. The following set:
\begin{equation}
\label{pf2_36}
\mathrm{Ext}\ K := K \cup \left(
\cup_{j=1}^n (K + \vec e_j)
\right),
\end{equation}
is said to be \textbf{the close extension of $K$}.
\end{definition}
The latter notion is obviously related to the notion of the closure of a set of monomials in $n$ variables
from~\cite[p. 89]{cit_988_Laurent_Mourrain}.
We think that this is a natural way to define extensions of prescribed moments in the case of general truncations.
Notice that the close extension of $K$ is also admissible.
Of course, in the case of triangular truncations we have
$$ \mathrm{Ext}\ K_r = K_{r+1},\qquad r\in\mathbb{Z}_+. $$

In contrast with Definition~\ref{d2_1} the following notion is not directly related to the notion of a border 
of a set of monomials in~\cite[p. 89]{cit_988_Laurent_Mourrain}.
\begin{definition}
\label{d2_2}
Let $K\subset\mathbb{Z}^n_+$ be an arbitrary admissible finite set; $n\in\mathbb{N}$. 
The following set:
\begin{equation}
\label{pf2_36_4}
\partial K := \{ \mathbf{m}\in K:\ \exists j\in \{1, ...,n\} \mbox{ such that } \mathbf{m} + \vec e_j\notin K \},
\end{equation}
is said to be \textbf{the border of $K$}.
\end{definition}
In the case of triangular truncations we have
$$ \partial K_r = K_{r}\backslash K_{r-1},\qquad r\in\mathbb{N}. $$
Choose $n=2$, $K=K_3\backslash\{ (1,1) \}$. Then we have $\mathrm{Ext}\ K = K_4$, and
$$ (\mathrm{Ext}\ K)\backslash K = \{ (1,1), (0,4), (1,3), (2,2), (3,1), (4,0) \} \nsubseteq \partial (\mathrm{Ext}\ K ) = $$
$$ = \{ (0,4), (1,3), (2,2), (3,1), (4,0) \}. $$
On the other hand, for triangular truncations we have
$$ (\mathrm{Ext}\ K_r)\backslash K_r = \partial (\mathrm{Ext}\ K_r ),\qquad r\in\mathbb{Z}_+. $$

\begin{theorem}
\label{t2_2}
Let $n\in\mathbb{N}$, $K\subset\mathbb{Z}^n_+$ be an arbitrary admissible finite set such that
\begin{equation}
\label{pf2_36_6}
(\mathrm{Ext}\ K)\backslash K \subseteq \partial (\mathrm{Ext}\ K ).
\end{equation}
Set $\mathcal{K} = K+K$, and $\mathcal{K}' = \mathrm{Ext}\ K + \mathrm{Ext}\ K$.
Let the moment problem~(\ref{f1_1}) with some prescribed moments $\mathcal{S} = (s_{\mathbf{k}})_{\mathbf{k}\in\mathcal{K}}$ be given,
and conditions~(\ref{f2_25}),(\ref{f2_45}) be satisfied (with some indexation of $K$).
Fix a continuation of the indexation of $K$ to $\mathrm{Ext}\ K$.
There exists a one-to-one correspondence between a set of all canonical solutions to the moment problem~(\ref{f1_1})
and a set of all dimensionally stable sets $\widehat{ \mathcal{S} } = 
(\widehat s_{\mathbf{k}})_{\mathbf{k}\in \mathcal{K}' }$, where 
$\widehat s_{\mathbf{k}}\in\mathbb{R}$, and
$\widehat s_{\mathbf{m}} = s_{\mathbf{m}}$, for
$\mathbf{m}\in \mathcal{K}$.
Here it is assumed that conditions~(\ref{f2_25}),(\ref{f2_45}) are satisfied for each $\widehat{ \mathcal{S} }$, with
the fixed indexation of $\mathrm{Ext}\ K$.

\end{theorem}
\textbf{Proof.}
Let $K$ have an indexation of type~(\ref{f2_9}). We extend this indexation to the whole set $\widehat K := \mathrm{Ext}\ K$:
\begin{equation}
\label{pf2_38}
\widehat K = \left\{
\mathbf{k}_j(\in\mathbb{Z}^n_+),\quad j=0,1,\ldots,\widehat\rho 
\right\}.
\end{equation}
The above indexation will be fixed throughout the whole proof.

Suppose that we have a dimensionally stable set $\widehat{ \mathcal{S} } = 
(\widehat s_{\mathbf{k}})_{\mathbf{k}\in \mathcal{K}' }$, where $\widehat s_{\mathbf{m}} = s_{\mathbf{m}}$, for
$\mathbf{m}\in \mathcal{K}$ (and conditions~(\ref{f2_25}),(\ref{f2_45}) are satisfied with
the fixed indexation of $\widehat K$).
Denote by $H$ and $\widehat H$ the associated Hilbert spaces for $\mathcal{S}$ and $\widehat{ \mathcal{S} }$, respectively.
In a similar manner, denote $\mathfrak{L}$, $\Omega_0$, $H_0$, $M_k$ ($k=1,...,n$) for $\mathcal{S}$ as in the Introduction, 
and the corresponding objects for $\widehat{ \mathcal{S} }$ denote by $\widehat{ \mathfrak{L} }$, $\widehat\Omega_0$, $\widehat H_0$, $\widehat M_k$.
Since $\widehat{ \mathcal{S} }$ is dimensionally stable, we have $\widehat H_0 = \widehat H$.
By Theorem~\ref{t2_1} we conclude that operators $\widehat M_k$, $k=1,...,n$, are commuting self-adjoint operators
acting in $\widehat H (=\widehat H_0)$.
Consider the following transformation $W$ from $H$ to $\widehat H$:
\begin{equation}
\label{pf2_42}
W \left[
\sum_{j=0}^\rho \xi_j \mathbf{t}^{ \mathbf{k}_j }
\right]_\mathfrak{L}
=
\left[
\sum_{j=0}^\rho \xi_j \mathbf{t}^{ \mathbf{k}_j }
\right]_{ \widehat{ \mathfrak{L} } },\qquad \xi_j\in\mathbb{C}.
\end{equation}
The transformation $W$ is well-defined and isometric.
In fact, suppose that
$$ h = 
\left[
\sum_{j=0}^\rho \xi_j \mathbf{t}^{ \mathbf{k}_j }
\right]_\mathfrak{L}
=
\left[
\sum_{j=0}^\rho \eta_j \mathbf{t}^{ \mathbf{k}_j }
\right]_\mathfrak{L},\quad \xi_j,\eta_j\in\mathbb{C}.
$$
Then
$$ \left\| 
\left[
\sum_{j=0}^\rho \xi_j \mathbf{t}^{ \mathbf{k}_j }
\right]_{ \widehat{ \mathfrak{L} } }
-
\left[
\sum_{j=0}^\rho \eta_j \mathbf{t}^{ \mathbf{k}_j }
\right]_{ \widehat{ \mathfrak{L} } }
\right\|^2 = $$
$$ = \sum_{j,l=0}^\rho (\xi_j -\eta_j) \overline{ (\xi_l -\eta_l) } 
\left(
\left[
\mathbf{t}^{ \mathbf{k}_j }
\right]_{ \widehat{ \mathfrak{L} } },
\left[
\mathbf{t}^{ \mathbf{k}_l }
\right]_{ \widehat{ \mathfrak{L} } }
\right)_{\widehat H} = 
\sum_{j,l=0}^\rho (\xi_j -\eta_j) \overline{ (\xi_l -\eta_l) } 
s_{ \mathbf{k}_j + \mathbf{k}_l } = $$
$$ = \sum_{j,l=0}^\rho (\xi_j -\eta_j) \overline{ (\xi_l -\eta_l) } 
\left(
\left[
\mathbf{t}^{ \mathbf{k}_j }
\right]_{ \mathfrak{L}  },
\left[
\mathbf{t}^{ \mathbf{k}_l }
\right]_{ \mathfrak{L} }
\right)_{H} = \| h-h \| = 0. $$
If
$$ g = 
\left[
\sum_{l=0}^\rho \gamma_l \mathbf{t}^{ \mathbf{k}_l }
\right]_\mathfrak{L},\qquad \gamma_l\in\mathbb{C},
$$
then
$$ (Wh, Wg) = \sum_{j,l=0}^\rho \xi_j \overline{ \gamma_l }
\left(
\left[
\mathbf{t}^{ \mathbf{k}_j }
\right]_{ \widehat{ \mathfrak{L} } },
\left[
\mathbf{t}^{ \mathbf{k}_l }
\right]_{ \widehat{ \mathfrak{L} } }
\right)_{\widehat H} = $$
$$ = \sum_{j,l=0}^\rho \xi_j \overline{ \gamma_l }
s_{ \mathbf{k}_j + \mathbf{k}_l } = 
\sum_{j,l=0}^\rho \xi_j \overline{ \gamma_l }
\left(
\left[
\mathbf{t}^{ \mathbf{k}_j }
\right]_{ \mathfrak{L} },
\left[
\mathbf{t}^{ \mathbf{k}_l }
\right]_{ \mathfrak{L} }
\right)_{H} = (h,g). $$

By the definition of the close extension of $K$ we obtain that $\{ 0,1,...,\rho \} \subseteq \widehat\Omega_0$. Therefore
$$ WH \subseteq \widehat H_0. $$
Notice that the indices $\{ \rho+1,...,\widehat\rho \}$ correspond to the elements of 
$(\mathrm{Ext}\ K)\backslash K$. By condition~(\ref{pf2_36_6}) these elements belong to $\partial (\mathrm{Ext}\ K )$.
Therefore they do not belong to $\widehat \Omega_0$, by the definition of the border.
Thus, we have $\{ 0,1,...,\rho \} = \widehat\Omega_0$, and
\begin{equation}
\label{pf2_44}
WH = \widehat H_0.
\end{equation}
The following operators are commuting self-adjoint operators in $H$:
\begin{equation}
\label{pf2_46}
C_k = W^{-1} \widehat M_k W,\qquad k=1,...,n.
\end{equation}
They generate a solution $\mu$ of the moment problem as it was described in the proof
of Proposition~1 in~\cite{cit_20000_Zagorodnyuk_2019_Concr_Oper}. Of course, this solution is canonical.

For each dimensionally stable set $\widehat{ \mathcal{S} }$ we put into correspondence a canonical solution $\mu$, which 
is constructed in the above manner.
We shall denote this map by~$\tau$: 
\begin{equation}
\label{pf2_48}
\tau \widehat{ \mathcal{S} } = \mu. 
\end{equation}
We are going to check that $\tau$ is a bijection between a set of all dimensionally stable close extensions of $\mathcal{S}$ and
a set of all canonical solutions of the moment problem.
At first, let us check that $\tau$ is surjective.

Let $\mu_c$ be an arbitrary canonical solution of the moment problem~(\ref{f1_1}).
There exist commuting self-adjoint operators $A_k\supseteq M_k$, $k=1,...,n$, in the associated Hilbert space $H$, and
\begin{equation}
\label{pf2_50}
\mu_c(\delta) = \left(
E(\delta) [1]_{\mathfrak{L}}, [1]_{\mathfrak{L}}
\right)_H,\qquad \delta\in\mathfrak{B}(\mathbb{R}^n),
\end{equation}
where $E(\delta)$ is the spectral measure of a commuting tuple $A_1,\ldots,A_n$ (see formula~(24) in~\cite{cit_20000_Zagorodnyuk_2019_Concr_Oper}).
Consider a set $\widetilde{ \mathcal{S} } := 
(\widetilde s_{\mathbf{k}})_{\mathbf{k}\in \mathcal{K}' }$, where 
\begin{equation}
\label{pf2_52}
\widetilde s_{\mathbf{k}} = \int \mathbf{t}^{\mathbf{k}} d\mu_c,\qquad \mathbf{k}\in \mathcal{K}'.
\end{equation}
Of course, $\widetilde s_{\mathbf{m}} = s_{\mathbf{m}}$, for
$\mathbf{m}\in \mathcal{K}$.
We claim that: $(1)$  $\widetilde{ \mathcal{S} }$ is dimensionally stable; $(2)$ $\tau \widetilde{ \mathcal{S} } = \mu_c$.
At first we observe that conditions~(\ref{f2_25}),(\ref{f2_45}) are satisfied for $\widetilde{ \mathcal{S} }$, with
the fixed indexation of $\mathrm{Ext}\ K$, since the extended moment problem~(\ref{pf2_52}) is solvable and the conditions are necessary
for the solvability. 
Let $\widetilde H$ be the associated Hilbert space for $\widetilde{ \mathcal{S} }$. 
The corresponding objects 
from the Introduction for $\widetilde{ \mathcal{S} }$ we denote by 
$\widetilde{ \mathfrak{L} }$, $\widetilde\Omega_0$, $\widetilde H_0$, $\widetilde M_k$.
Notice that (see formula~(22) in~\cite{cit_20000_Zagorodnyuk_2019_Concr_Oper})
\begin{equation}
\label{pf2_54}
\mathop{\rm Lin}\nolimits \left\{ A_1^{k_1} A_2^{k_2} ... A_n^{k_n} [1]_\mathfrak{L},\ (k_1,...,k_n)\in\mathbb{Z}_+ \right\} = H.
\end{equation}
Denote by $L^2_{\mu_c}$ the space of all (classes of the equivalence) of all square integrable complex-valued $\mathfrak{B}(\mathbb{R}^n)$-measurable 
functions $f(\mathbf{t})$, $\mathbf{t}\in\mathbb{R}^n$. 
The classes of the equivalence we shall denote by $[f]=[f]_{ L^2_{\mu_c} }$. Set
\begin{equation}
\label{pf2_56}
L^2_{\mu_c,0} = 
\mathop{\rm Lin}\nolimits \left\{ \left[ t_1^{k_1} t_2^{k_2} ... t_n^{k_n} \right]_{ L^2_{\mu_c} },\ (k_1,...,k_n)\in\mathbb{Z}^n_+ \right\}.
\end{equation}
Consider the following transformation $V$ acting from $H$ to $L^2_{\mu_c,0}$:
\begin{equation}
\label{pf2_58}
V h =
\sum_{ \mathbf{k}\in\mathbb{Z}_+ } \alpha_{\mathbf{k}}
\left[ t_1^{k_1} t_2^{k_2} ... t_n^{k_n} \right]_{ L^2_{\mu_c} },\quad 
h = \sum_{ \mathbf{k}\in\mathbb{Z}_+ } \alpha_{\mathbf{k}} A_1^{k_1} A_2^{k_2} ... A_n^{k_n} [1]_\mathfrak{L},\
\alpha_{\mathbf{k}}\in\mathbb{C},
\end{equation}
where all but finite number of $\alpha_{\mathbf{k}}$s are zero and this will be assumed in similar situations in what follows.
Let us check that $V$ is well-defined. Suppose that an element $h\in H$ has another representation:
\begin{equation}
\label{pf2_60}
h = \sum_{ \mathbf{k}\in\mathbb{Z}_+ } \beta_{\mathbf{k}} A_1^{k_1} A_2^{k_2} ... A_n^{k_n} [1]_\mathfrak{L},\
\beta_{\mathbf{k}}\in\mathbb{C}.
\end{equation}
Using relation~(\ref{pf2_50}) we may write:
$$ \left\|
\sum_{ \mathbf{k}\in\mathbb{Z}_+ } (\alpha_{\mathbf{k}} - \beta_{\mathbf{k}})
\left[ t_1^{k_1} t_2^{k_2} ... t_n^{k_n} \right]_{ L^2_{\mu_c} }
\right\|^2 = $$
$$ = \sum_{ \mathbf{k},\mathbf{m} \in\mathbb{Z}_+ } (\alpha_{\mathbf{k}} - \beta_{\mathbf{k}})
\overline{ (\alpha_{\mathbf{m}} - \beta_{\mathbf{m}}) }
\left(
\left[ t_1^{k_1} t_2^{k_2} ... t_n^{k_n} \right]_{ L^2_{\mu_c} },
\left[ t_1^{m_1} t_2^{m_2} ... t_n^{m_n} \right]_{ L^2_{\mu_c} }
\right)
 = $$
$$ = \sum_{ \mathbf{k},\mathbf{m} \in\mathbb{Z}_+ } (\alpha_{\mathbf{k}} - \beta_{\mathbf{k}})
\overline{ (\alpha_{\mathbf{m}} - \beta_{\mathbf{m}}) }
\int \mathbf{t}^{ \mathbf{k}+\mathbf{m} } d\mu_c = $$
$$ = \sum_{ \mathbf{k},\mathbf{m} \in\mathbb{Z}_+ } (\alpha_{\mathbf{k}} - \beta_{\mathbf{k}})
\overline{ (\alpha_{\mathbf{m}} - \beta_{\mathbf{m}}) }
\left(
A_1^{k_1} A_2^{k_2} ... A_n^{k_n} [1]_\mathfrak{L},
A_1^{m_1} A_2^{m_2} ... A_n^{m_n} [1]_\mathfrak{L}
\right)_H
= $$
$$ = \left\|
\sum_{ \mathbf{k}\in\mathbb{Z}_+ } (\alpha_{\mathbf{k}} - \beta_{\mathbf{k}})
A_1^{k_1} A_2^{k_2} ... A_n^{k_n} [1]_\mathfrak{L}
\right\|^2 = \| h - h \|_H = 0. $$
Therefore $V$ is well-defined. If
\begin{equation}
\label{pf2_62}
g = \sum_{ \mathbf{m}\in\mathbb{Z}_+ } \gamma_{\mathbf{m}} A_1^{m_1} A_2^{m_2} ... A_n^{m_n} [1]_\mathfrak{L},\
\gamma_{\mathbf{m}}\in\mathbb{C},
\end{equation}
then
$$ (Vh,Vg) = 
\sum_{ \mathbf{k},\mathbf{m} \in\mathbb{Z}_+ } \alpha_{\mathbf{k}} \overline{ \gamma_{\mathbf{m}} }
\left(
\left[ t_1^{k_1} t_2^{k_2} ... t_n^{k_n} \right]_{ L^2_{\mu_c} },
\left[ t_1^{m_1} t_2^{m_2} ... t_n^{m_n} \right]_{ L^2_{\mu_c} }
\right) = 
$$
$$ = \sum_{ \mathbf{k},\mathbf{m} \in\mathbb{Z}_+ } \alpha_{\mathbf{k}} \overline{ \gamma_{\mathbf{m}} }
\int \mathbf{t}^{ \mathbf{k}+\mathbf{m} } d\mu_c = $$
$$ = \sum_{ \mathbf{k},\mathbf{m} \in\mathbb{Z}_+ } \alpha_{\mathbf{k}} \overline{ \gamma_{\mathbf{m}} }
\left(
A_1^{k_1} A_2^{k_2} ... A_n^{k_n} [1]_\mathfrak{L},
A_1^{m_1} A_2^{m_2} ... A_n^{m_n} [1]_\mathfrak{L}
\right) = (h,g). $$
Thus, $V$ is a linear isometric transformation which maps $H$ onto $L^2_{\mu_c,0}$. Then
\begin{equation}
\label{pf2_64}
\mathrm{dim}\ L^2_{\mu_c,0} = \mathrm{dim}\ H.
\end{equation}

Consider the following transformation $\widetilde W$ from $H$ to $\widetilde H$:
\begin{equation}
\label{pf2_62_4}
\widetilde W \left[
\sum_{j=0}^\rho \xi_j \mathbf{t}^{ \mathbf{k}_j }
\right]_\mathfrak{L}
=
\left[
\sum_{j=0}^\rho \xi_j \mathbf{t}^{ \mathbf{k}_j }
\right]_{ \widetilde{ \mathfrak{L} } },\qquad \xi_j\in\mathbb{C}.
\end{equation}
The transformation $\widetilde W$ is well-defined and isometric.
This can be checked in the same manner, as it was done for the transformation $W$.
Moreover, we also have an analogue of~(\ref{pf2_44}):
\begin{equation}
\label{pf2_62_5}
\widetilde W H = \widetilde H_0.
\end{equation}
Then
\begin{equation}
\label{pf2_62_7}
\mathrm{dim}\ H = \mathrm{dim}\ \widetilde W H = \mathrm{dim}\ \widetilde H_0 \leq \mathrm{dim}\ \widetilde H.
\end{equation}

Now consider the following transformation $U$, acting from $\widetilde H$ into $L^2_{\mu,0}$:
\begin{equation}
\label{pf2_70}
U h =
\sum_{j=0}^{ \widehat{\rho} } \alpha_j
\left[ {\mathbf{t}}^{ \mathbf{k}_j } \right]_{ L^2_{\mu_c} },\quad 
h = \sum_{j=0}^{ \widehat{\rho} } \alpha_j \left[ {\mathbf{t}}^{ \mathbf{k}_j }
\right]_{ \widetilde{ \mathfrak{L} } },\
\alpha_j\in\mathbb{C}.
\end{equation}
If $h$ has another representation:
$$ h = \sum_{j=0}^{ \widehat{\rho} } \beta_j \left[ {\mathbf{t}}^{ \mathbf{k}_j } \right]_{ \widetilde{ \mathfrak{L} } },\quad
\beta_j\in\mathbb{C}, $$
then
$$ \left\| 
\sum_{j=0}^{ \widehat{\rho} } (\alpha_j - \beta_j)
\left[ {\mathbf{t}}^{ \mathbf{k}_j } \right]_{ L^2_{\mu_c} }
\right\|^2 = \sum_{j,m=0}^{ \widehat{\rho} } (\alpha_j - \beta_j) \overline{ (\alpha_m - \beta_m) }
\left(
\left[ {\mathbf{t}}^{ \mathbf{k}_j } \right]_{ L^2_{\mu_c} },
\left[ {\mathbf{t}}^{ \mathbf{k}_m } \right]_{ L^2_{\mu_c} }
\right) = $$
$$ = \sum_{j,m=0}^{ \widehat{\rho} } (\alpha_j - \beta_j) \overline{ (\alpha_m - \beta_m) }
\int  {\mathbf{t}}^{ \mathbf{k}_j + \mathbf{k}_m } d\mu_c = 
\sum_{j,m=0}^{ \widehat{\rho} } (\alpha_j - \beta_j) \overline{ (\alpha_m - \beta_m) }
\widetilde s_{ \mathbf{k}_j + \mathbf{k}_m }  = $$
$$ = \sum_{j,m=0}^{ \widehat{\rho} } (\alpha_j - \beta_j) \overline{ (\alpha_m - \beta_m) }
\left(
\left[ {\mathbf{t}}^{ \mathbf{k}_j } \right]_{ \widetilde{ \mathfrak{L} } },
\left[ {\mathbf{t}}^{ \mathbf{k}_m } \right]_{ \widetilde{ \mathfrak{L} } }
\right)_{ \widetilde H }
= \| h - h \|^2 = 0 . $$
Thus, $U$ is well-defined. If 
$$ g = \sum_{j=0}^{ \widehat{\rho} } \gamma_j \left[ {\mathbf{t}}^{ \mathbf{k}_j } \right]_{ \widetilde{ \mathfrak{L} } },\
\gamma_j\in\mathbb{C}, $$
then
$$ (Uh, Ug) = \sum_{j,m=0}^{ \widehat{\rho} } \alpha_j \overline{ \gamma_m }
\left(
\left[ {\mathbf{t}}^{ \mathbf{k}_j } \right]_{ L^2_{\mu_c} },
\left[ {\mathbf{t}}^{ \mathbf{k}_m } \right]_{ L^2_{\mu_c} }
\right) = $$
$$ = \sum_{j,m=0}^{ \widehat{\rho} } \alpha_j \overline{ \gamma_m }
\widetilde s_{ \mathbf{k}_j + \mathbf{k}_m } = 
\sum_{j,m=0}^{ \widehat{\rho} } \alpha_j \overline{ \gamma_m }
\left(
\left[ {\mathbf{t}}^{ \mathbf{k}_j } \right]_{ \widetilde{ \mathfrak{L} } },
\left[ {\mathbf{t}}^{ \mathbf{k}_m } \right]_{ \widetilde{ \mathfrak{L} } }
\right)_{ \widetilde H }
= (h,g). $$
So, $U$ is a linear isometric transformation which maps $\widetilde H$ onto a subspace of
$L^2_{\mu_c,0}$. It follows that
\begin{equation}
\label{pf2_72}
\mathrm{dim}\ \widetilde H \leq \mathrm{dim}\ L^2_{\mu_c,0}.
\end{equation}
By~(\ref{pf2_64}),(\ref{pf2_62_7}),(\ref{pf2_72}) we obtain that
\begin{equation}
\label{pf2_74}
\mathrm{dim}\ \widetilde H =  \mathrm{dim}\ \widetilde H_0 = \mathrm{dim}\ L^2_{\mu_c,0} = \mathrm{dim}\ H.
\end{equation}
So, $\widetilde{ \mathcal{S} }$ is dimensionally stable and the first claim is proved.
Notice that relation~(\ref{pf2_74}) also shows that $U$ is a unitary transformation which maps $\widetilde H$ onto $L^2_{\mu_c,0}$.

\noindent
Now we are going to check the second claim. Since $\widetilde{ \mathcal{S} }$ is dimensionally stable,
we can construct a canonical solution $\widetilde\mu := \tau \widetilde{ \mathcal{S} }$. 
Using the same notations as at the beginning of the proof, we can repeat all the constructions for the concrete choice of moments
$\widehat{ \mathcal{S} } = \widetilde{ \mathcal{S} }$. The use of the same notations will cause no confusion,
since we shall not use these constructions elsewhere. So, we now have double notations for the
same objects: $\widehat\Omega_0 = \widetilde\Omega_0$, $\widehat M_k = \widetilde M_k$, etc. 
We have
\begin{equation}
\label{pf2_76}
\widetilde \mu(\delta) = \left(
E_C(\delta) [1]_{ \mathfrak{L} }, [1]_{ \mathfrak{L} }
\right)_H,\qquad \delta\in\mathfrak{B}(\mathbb{R}^n),
\end{equation}
where $E_C(\delta)$ is the spectral measure of the commuting tuple $C_1,\ldots, C_n$ ($C_k = W^{-1} \widetilde M_k W$).
In order to prove the second claim we need to check that $\widetilde\mu = \mu_c$.
Choose an arbitrary element $h\in H$ of the form~(\ref{pf2_60}). We may write
$$ V A_l h = 
V \sum_{ \mathbf{k}\in\mathbb{Z}_+ } \beta_{\mathbf{k}} A_1^{k_1}\cdots A_l^{k_l+1}\cdots A_n^{k_n} [1]_\mathfrak{L} =
$$
$$ = \sum_{ \mathbf{k}\in\mathbb{Z}_+ } \beta_{\mathbf{k}}
\left[ t_1^{k_1} \cdots  t_l^{k_l+1} \cdots  t_n^{k_n} \right]_{ L^2_{\mu_c} } = \Lambda_l V h, $$
where we denote by $\Lambda_l$ the operator of multiplication by $t_l$ in $L^2_{\mu_c}$, restricted to $L^2_{\mu_c,0}$.
Therefore
\begin{equation}
\label{pf2_80}
V A_l V^{-1} = \Lambda_l,\qquad l=1,...,n.
\end{equation}
Choose an arbitrary element $g\in \widetilde H_0 (=\widetilde H)$, 
$$ g = \sum_{j\in\widetilde\Omega_0} \alpha_j \left[ \mathbf{t}^{\mathbf{k}_j} \right]_{\widetilde{\mathfrak{L}}},\qquad \alpha_j\in\mathbb{C}. $$
Then
$$ U \widetilde M_l g = U \sum_{j\in\widetilde\Omega_0} \alpha_j \left[ \mathbf{t}^{ \mathbf{k}_j + \vec e_l } \right]_{\widetilde{\mathfrak{L}}} = 
\sum_{j\in\widetilde\Omega_0} \alpha_j \left[ 
\mathbf{t}^{ \mathbf{k}_j + \vec e_l }
\right]_{L^2_{\mu_c}} = \Lambda_l U g.
$$
Therefore
\begin{equation}
\label{pf2_82}
U \widetilde M_l U^{-1} = \Lambda_l,\qquad l=1,...,n.
\end{equation}
By~(\ref{pf2_80}),(\ref{pf2_82}) we obtain that $A_l$ and $\widetilde M_l$ are unitary equivalent:
\begin{equation}
\label{pf2_84}
T A_l T^{-1} = \widetilde M_l,\qquad l=1,...,n,
\end{equation}
where $T := U^{-1} V$.
Therefore operators $A_l$ and $C_l$ are also unitarily equivalent:
\begin{equation}
\label{pf2_84_3}
F A_l F^{-1} = C_l,\qquad l=1,...,n,
\end{equation}
where $F := W^{-1} T$.
Operators $A_l$ and $C_l$ have the same eigenvalues which we denote by
$$ \{ \lambda_{l;j} \}_{j=1}^{n_l},\quad \mbox{where $\lambda_{l;m}\not= \lambda_{l;k}$, if $m\not= k$.} $$
For $\lambda_{l;j}$ we denote by $H_{l;j}$ and $\widetilde H_{l;j}$ the eigen subspaces of $A_l$ and $C_l$, respectively. 
We also denote $P_{l,j} := P^H_{H_{l;j}}$ and $\widetilde P_{l,j} := P^H_{\widetilde H_{l;j}}$.
Then
\begin{equation}
\label{pf2_86}
\widetilde P_{l,j} = F P_{l,j} F^{-1},\qquad j\in\mathbb{Z}_{1,n_l},\ l\in\mathbb{Z}_{1,n}.
\end{equation}
The measures $\mu_c$ and $\widetilde\mu$ are atomic and they can have atoms at the following points only:
\begin{equation}
\label{pf2_88}
\mathbf{a}_{j_1,j_2,...,j_n} := ( \lambda_{1;j_1}, \lambda_{2;j_2},...,\lambda_{n;j_n} ),\qquad  j_k\in\mathbb{Z}_{1,n_k},\ k\in\mathbb{Z}_{1,n}.
\end{equation}
The corresponding masses have the following values: 
$$ \mu_c (\{ \mathbf{a}_{j_1,j_2,...,j_n} \}) = 
\left(
E_1 (\{ \lambda_{1;j_1} \} ) \cdots E_n (\{ \lambda_{n;j_n} \} ) [1]_{\mathfrak{L}},
[1]_{\mathfrak{L}}
\right)_H 
= $$
\begin{equation}
\label{pf2_90}
= \left(
P_{1,j_1} \cdots P_{n,j_n} [1]_{\mathfrak{L}}, [1]_{\mathfrak{L}}
\right)_H,
\end{equation}
$$ \widetilde \mu (\{ \mathbf{a}_{j_1,j_2,...,j_n} \}) = 
\left(
E_{C_1} (\{ \lambda_{1;j_1} \} ) \cdots E_{C_n} (\{ \lambda_{n;j_n} \} ) [1]_{\mathfrak{L}},
[1]_{\mathfrak{L}}
\right)_H 
= $$
\begin{equation}
\label{pf2_92}
= \left(
\widetilde P_{1,j_1} \cdots \widetilde P_{n,j_n} [1]_{\mathfrak{L}},
[1]_{\mathfrak{L}}
\right)_H,
\end{equation}
where $E_j$ and $E_{C_j}$ are the spectral measures of $A_j$ and $C_j$, respectively.
By~(\ref{pf2_86}) we conclude that $\mu_c = \widetilde\mu$. Thus, the second claim is proved.
We conclude that the map $\tau$ is surjective. Our considerations also show that the set of canonical solutions is empty
if and only if the set of dimensionally stable extensions is empty.

Suppose that we have two different dimensionally stable sets
$\mathcal{S}' = (s_{\mathbf{k}}')_{\mathbf{k}\in \mathcal{K}' }$, 
and
$\mathcal{S}'' = (s_{\mathbf{k}}'')_{\mathbf{k}\in \mathcal{K}' }$, 
where $s_{\mathbf{m}}' = s_{\mathbf{m}}'' = s_{\mathbf{m}}$, for $\mathbf{m}\in \mathcal{K}$.
Denote
$$ \mu_1 := \tau \mathcal{S}',\quad \mu_2 := \tau \mathcal{S}''. $$
Since $\mathcal{S}'\not = \mathcal{S}''$, then there exists $\mathbf{k}\in \mathcal{K}'$ such that
$s_{\mathbf{k}}'\not = s_{\mathbf{k}}''$.
Then
$$ \int \mathbf{t}^{\mathbf{k}} d\mu_1 \not= \int \mathbf{t}^{\mathbf{k}} d\mu_2. $$
Therefore $\mu_1\not= \mu_2$.
Thus, $\tau$ is injective.     
It follows that $\tau$ is a biection and this completes the proof of the theorem.
$\Box$

\section{The two-dimensional case with triangular truncations. Indices of nonself-adjointness.}

Fix an arbitrary number $r\in\mathbb{N}$ and consider $K = K_r$, where $K_r$ is given by~(\ref{pf2_7_5}) with $n=2$.
Assume that the moment problem~(\ref{f1_1}) with 
$\mathcal{K} = K+K$ and some $\mathcal{S} = (s_{\mathbf{k}})_{\mathbf{k}\in\mathcal{K}}$ is given and $s_{ (0,0) } > 0$
(this condition excludes the trivial case).
Choose and fix an indexation~(\ref{f2_9}) such that the elements of $K_{r-1}$ are indexed at first. 
Assume that conditions~(\ref{f2_25}),(\ref{f2_45})
are satisfied.
Construct the associated Hilbert space $H$, symmetric operators $M_1$, $M_2$, and other related objects from the Introduction.
Notice that the operators $M_1$, $M_2$ are now defined on the same subspace $H_0$.
As it was already mentioned in the Introduction, the number
\begin{equation}
\label{pf3_5}
i_s := \mathrm{dim}\ (H\ominus H_0),
\end{equation}
is said to be \textit{the index of nonself-adjointness of $\mathcal{S}$}.
The case $i_s=0$ leads to the dimension stability or flatness of $\mathcal{S}$. 
From the definitions of $K=K_r$ and $\Omega_0$ it is clear that $i_s$ can take values
$0, 1, ..., r+1$.

In what follows we shall assume that $i_s > 0$. 
Apply the Gram-Schmidt orthogonalization procedure to the sequence $\left[ \mathbf{t}^{ \mathbf{k}_j } \right]_{\mathfrak{L}}$, $j=0,1,...,\rho$.
We shall obtain an orthonormal basis $\mathcal{G} = \{ g_j \}_{j=0}^{d+i_s}$ in $H$, with some $d\in\mathbb{Z}_+$. Of course, $\mathrm{dim}\ H_0 = d+1$,
$\mathrm{dim}\ H = d+i_s+1$, and $\{ g_j \}_{j=0}^d$ is an orthonormal basis in $H_0$.
Operators $M_1$ and $M_2$ are defined on $H_0$, which is \textit{a proper subspace of $H$}.
The existence of self-adjoint extensions of such symmetric operators was investigated by Krasnoselskii, for details see, e.g.,
a survey in~\cite{cit_14550_Zagorodnyuk_2013_Survey__Generalized resolvents} and references therein.
Of course, the existence of self-adjoint extensions in a finite-dimensional case is quite clear from the matrix representations
of the corresponding operators.

Let $R_1$ and $R_2$, $D(R_1) = D(R_2) = H$, be arbitrary self-adjoint extensions of $M_1$ and $M_2$, respectively. 
For the basis $\mathcal{G}$ we denote by
$\mathcal{R}_1$ and $\mathcal{R}_2$ the matrices of $R_1$ and $R_2$, respectively. We have
\begin{equation}
\label{pf3_7}
\mathcal{R}_k = \left(
(R_k g_l, g_j)_H 
\right)_{j,l=0}^{ d+i_s } =
\left( 
\begin{array}{cc} A_k & B_k^* \\
B_k & C_k \end{array}
\right),\quad k=1,2,
\end{equation}
where 
\begin{equation}
\label{pf3_9}
A_k = \left(
(R_k g_l, g_j)_H 
\right)_{j,l=0}^d =
\left(
(M_k g_l, g_j)_H 
\right)_{j,l=0}^d,
\end{equation}
\begin{equation}
\label{pf3_12}
B_k = \left(
(R_k g_l, g_j)_H 
\right)_{j=d+1,...,d+i_s;\ l=0,...,d}
=
\left(
(M_k g_l, g_j)_H 
\right)_{j=d+1,...,d+i_s;\ l=0,...,d},
\end{equation}
\begin{equation}
\label{pf3_14}
C_k = \left(
(R_k g_l, g_j)_H 
\right)_{j,l=d+1}^{d+i_s}.
\end{equation}
Thus, matrices $A_k$ and $B_k$ are calculated by the given moments.
A direct block multiplication of the corresponding matrices shows that the operators $R_1$ and $R_2$ commute 
if and only if the following conditions hold:
\begin{equation}
\label{pf3_16}
A_1 A_2 + B_1^* B_2 = A_2 A_1 + B_2^* B_1,
\end{equation}
\begin{equation}
\label{pf3_18}
B_2^* C_1 - B_1^* C_2 = A_1 B_2^* - A_2 B_1^*,
\end{equation}
\begin{equation}
\label{pf3_20}
C_1 C_2 - C_2 C_1 = B_2 B_1^* - B_1 B_2^*.
\end{equation}
Observe that relation~(\ref{pf3_16}) is a necessary condition for the existence of canonical solutions.
By~(\ref{pf3_20}) we see that the following condition:
\begin{equation}
\label{pf3_22}
\mathop{\rm Tr}\nolimits (B_2 B_1^* - B_1 B_2^*) = 0,
\end{equation}
is necessary for the existence of canonical solutions as well.

\begin{proposition}
\label{p3_1}
Let $n=2$, $r\in\mathbb{N}$ and $K = K_r$, with $K_r$ as in~(\ref{pf2_7_5}).
Assume that the moment problem~(\ref{f1_1}) with $\mathcal{K} = K+K$
and some $\mathcal{S} = (s_{\mathbf{k}})_{\mathbf{k}\in\mathcal{K}}$ is given with $s_{ (0,0) } > 0$.
Choose and fix an indexation~(\ref{f2_9}) such that the elements of $K_{r-1}$ are indexed at first, an then the rest of $K_r$. 
Assume that conditions~(\ref{f2_25}),(\ref{f2_45})
are satisfied.
Construct the associated Hilbert space~$H$, the symmetric operators $M_1,M_2$, and define other related objects as in the Introduction. 
Suppose that $i_s > 0$ and conditions~(\ref{pf3_16}),(\ref{pf3_22}) hold with $A_k,B_k$, $k=1,2$, defined by~(\ref{pf3_9}),(\ref{pf3_12}).
The moment problem~(\ref{f1_1}) has a canonical solution if and only if
conditions~(\ref{pf3_18}),(\ref{pf3_20}) hold for some Hermitian matrices $C_1,C_2$ of size $i_s\times i_s$.
\end{proposition}
\textbf{Proof.} 
\textit{Necessity.} Suppose that there exists a canonical solution $\widetilde \mu$ of the moment problem~(\ref{f1_1}).
It is generated by commuting self-adjoint extensions $\widetilde R_j\supseteq M_j$, $j=1,2$.
Applying the considerations before the statement of the theorem for $R_j = \widetilde R_j$, we obtain that
conditions~(\ref{pf3_18}),(\ref{pf3_20}) hold for the corresponding Hermitian matrices $C_1,C_2$.

\noindent
\textit{Sufficiency.}
Assume that there exist Hermitian matrices $C_1,C_2$ of size $i_s\times i_s$ which satisfy~(\ref{pf3_18}),(\ref{pf3_20}).
Define operators $R_k$. $k=1,2$, by the matrices $\mathcal{R}_k$, as in~(\ref{pf3_7}) with
$$ A_k = \left(
(M_k g_l, g_j)_H 
\right)_{j,l=0}^d,\quad
B_k = \left(
(M_k g_l, g_j)_H 
\right)_{j=d+1,...,d+i_s;\ l=0,...,d}. $$
Operators $R_1$, $R_2$ are commuting self-adjoint operators, extending $M_1$ and $M_2$, respectively. 
Operators $R_k$ generate a canonical solution $\mu$ of the moment problem.
The proof is complete.
$\Box$

Suppose that assumptions of Proposition~\ref{p3_1} hold. Now we shall study equations~(\ref{pf3_18}),(\ref{pf3_20})
for unknown Hermitian matrices $C_1$, $C_2$ of size $i_s\times i_s$.
Let
\begin{equation}
\label{pf3_24}
C_k = (c_{k;j.l})_{j,l=1}^{i_s},\qquad c_{k;j.l}\in\mathbb{C},\quad k=1,2.
\end{equation}
For $j,l\in\mathbb{Z}_{1,i_s}:\ j>l$, we may write
\begin{equation}
\label{pf3_26}
c_{k;j.l} = \alpha_{k;j.l} + i \beta_{k;j.l},\qquad \alpha_{k;j.l},\beta_{k;j.l}\in\mathbb{R}. 
\end{equation}
Thus, Hermitian matrices $C_1$ and $C_2$ are determined by real numbers
\begin{equation}
\label{pf3_28}
\alpha_{k;j.l},\beta_{k;j.l},\qquad j,l\in\mathbb{Z}_{1,i_s}:\ j>l,
\end{equation}
and
\begin{equation}
\label{pf3_30}
c_{k;m.m},\qquad m\in\mathbb{Z}_{1,i_s};\qquad  k=1,2.
\end{equation}

In the case $i_s=1$  relation~(\ref{pf3_20}) takes the following form:
\begin{equation}
\label{pf3_32}
B_2 B_1^* - B_1 B_2^* = 0,
\end{equation}
and it can be verified directly.
In this case, taking the real and the imaginary parts of both sides of relation~(\ref{pf3_18}) we obtain
a system of linear algebraic equations with real coefficients, with respect to unknown real numbers $\alpha_{k;j.l},\beta_{k;j.l}$, and $c_{k;m.m}$.
Thus, in the case $i_s=1$ the existence of canonical solutions  can be easily checked.

Now assume that $i_s=2$. By~(\ref{pf3_22}) we may write:
\begin{equation}
\label{pf3_34}
B_2 B_1^* - B_1 B_2^* = U
\left(
\begin{array}{cc} ir & 0 \\
0 & -ir \end{array}
\right)
U^{-1},\qquad r\in\mathbb{R},
\end{equation}
where $U$ is a suitable unitary matrix.
Relation~(\ref{pf3_20}) may be written in the following form:
\begin{equation}
\label{pf3_36}
\mathbf{C}_1 \mathbf{C}_2 - \mathbf{C}_2 \mathbf{C}_1 = 
\left(
\begin{array}{cc} ir & 0 \\
0 & -ir \end{array}
\right),
\end{equation}
where
\begin{equation}
\label{pf3_38}
\mathbf{C}_1 = U^{-1} C_1 U,\quad  \mathbf{C}_2 = U^{-1} C_2 U.
\end{equation}
Relation~(\ref{pf3_18}) takes the following form:
\begin{equation}
\label{pf3_42}
\mathcal{B}_2 \mathbf{C}_1 - \mathcal{B}_1 \mathbf{C}_2 = \mathcal{D},
\end{equation}
where
\begin{equation}
\label{pf3_44}
\mathcal{B}_2 = B_2^* U,\quad \mathcal{B}_1 = B_1^* U,\quad
\mathcal{D} = ( A_1 B_2^* - A_2 B_1^* ) U.
\end{equation}
Thus, the moment problem~(\ref{f1_1}) has a canonical solution if and only if
conditions~(\ref{pf3_36}),(\ref{pf3_42}) hold for some Hermitian matrices $\mathbf{C}_1,\mathbf{C}_2$ of size $2\times 2$.
Let
\begin{equation}
\label{pf3_46}
\mathbf{C}_1 = 
\left(
\begin{array}{cc} a & c \\
\overline{c} & b \end{array}
\right),\quad 
\mathbf{C}_2 = 
\left(
\begin{array}{cc} d & g \\
\overline{g} & f \end{array}
\right),
\end{equation}
where $a,b,d,f\in\mathbb{R}$, $c=c'+i c''$, $g=g'+ig''$, $c',c'',g',g''\in\mathbb{R}$.
If we take the real and imaginary parts of relation~(\ref{pf3_42}), we shall obtain a system of linear algebraic equations
with real coefficients and real unknowns $a,b,d,f,c',c'',g',g''$. 
As for relation~(\ref{pf3_36}), it is equivalent to the following two equations:
\begin{equation}
\label{pf3_48}
c \overline{g} - \overline{c} g = i r,
\end{equation}
\begin{equation}
\label{pf3_50}
c (f-d) + g (a-b) = 0.
\end{equation}
Relations~(\ref{pf3_48}),(\ref{pf3_50}) are equivalent to the following three equations:
\begin{equation}
\label{pf3_52}
c''g' - c' g'' = \frac{r}{2},
\end{equation}
\begin{equation}
\label{pf3_54}
c' (f-d) + g' (a-b) = 0,
\end{equation}
\begin{equation}
\label{pf3_56}
c'' (f-d) + g'' (a-b) = 0.
\end{equation}

\noindent
\textbf{Case 1}: $r\not= 0$. In this case relations~(\ref{pf3_54}),(\ref{pf3_56}) imply $f=d$, $a=b$.
In fact, these relations form a real linear algebraic system of equations with respect to $f-d, a-b$,
having a non-zero determinant.
Thus we have a non-linear equation~(\ref{pf3_52}) and a linear system of equations including 
$f=d$, $a=b$, and those equations obtained from relation~(\ref{pf3_42}). 
Solvability of this system is necessary for the the existence of canonical solutions.

Suppose that the latter linear system of equations has solutions.
We substitute its solution into relation~(\ref{pf3_52}).
If the solution of the linear system was unique, we shall obtain a necessary condition for the existence of
canonical solutions.
Otherwise, we obtain a single quadratic equation with several unknowns. It can be effectively solved for the cases of one or two
unknowns. The case of three or more unknowns for a quadratic equation seems to be not investigated.

\noindent
\textbf{Case 2}: $r = 0$. In this case relations~(\ref{pf3_48}),(\ref{pf3_50}) take the following form:
\begin{equation}
\label{pf3_58}
c \overline{g} - \overline{c} g = 0,
\end{equation}
\begin{equation}
\label{pf3_60}
c (f-d) + g (a-b) = 0.
\end{equation}
Let us look for solutions of these equations with $g=0$.
Then $c(f-d) = 0$. Thus, in this special case equations~(\ref{pf3_58}),(\ref{pf3_60}) will be satisfied if $c=0$ or $f=d$.
It remains to solve the linear system corresponding to the real and imaginary parts of relation~(\ref{pf3_42}).

In a similar way we may look for solutions of equations~(\ref{pf3_58}),(\ref{pf3_60}) with $c=0$. 
In this special case equations~(\ref{pf3_58}),(\ref{pf3_60}) will be fulfilled if  $g=0$ or $a=b$.
Then we also come to a linear system of remaining equations.

Let us search for solutions of equations~(\ref{pf3_58}),(\ref{pf3_60}) with $g\not= 0$ and $c\not= 0$.
Rewrite equations~(\ref{pf3_58}),(\ref{pf3_60}) in the following form:
\begin{equation}
\label{pf3_62}
\frac{c}{g} \in\mathbb{R},
\end{equation}
\begin{equation}
\label{pf3_64}
\frac{c}{g} (f-d) + a-b = 0.
\end{equation}
These equations have a solution if and only if the following system with an additional real unknown $\beta$ is solvable:
\begin{equation}
\label{pf3_68}
c' = \beta g',\quad c'' = \beta g'',
\end{equation}
\begin{equation}
\label{pf3_70}
\beta (f-d) + a-b = 0.
\end{equation}
Together with the linear system, corresponding to the real and imaginary parts of relation~(\ref{pf3_42}),
we get a linear system with coefficients depending on a real parameter $\beta$.
In this case one can apply the Gauss elimination method. In fact, the choice of a leading element may lead to
a search of real roots of a polynomial in $\beta$. This roots can be localized with any desired precision.
If one gets some equations without unknowns, they also lead to a search of real roots of a polynomial in $\beta$.

\begin{corollary}
\label{c3_1}
Suppose that assumptions of Proposition~\ref{p3_1} hold.
If $i_s=1$ then the existence of canonical solutions can be checked explicitly.
If $i_s=2$ then the existence of canonical solutions either can be checked explicitly, or one gets a single quadratic
equation with several real unknowns.
\end{corollary}
\textbf{Proof.} The proof follows from considerations before the statement of the corollary.
$\Box$

Let us show that the cases $i_s=1$ and $i_s=2$ can really happen.

\noindent
\textbf{Example 1} ($i_s=1$).
Let $n=2$, $K = K_2$ (see~(\ref{pf2_7_5})), and $\mathcal{K} = K+K = K_4$.
Consider the moment problem~(\ref{f1_1}) with the moments $\mathcal{S} = (s_{\mathbf{k}})_{\mathbf{k}\in\mathcal{K}}$:

$$ s_{ (0,0) } = 9,\ s_{ (1,0) } = -1,\ s_{ (0,1) } = 0,\  s_{ (2,0) } = 1,\ s_{ (1,1) } = 0,\ s_{ (0,2) } = 2,\ $$

$$ s_{ (3,0) } = -1,\ s_{ (2,1) } = s_{ (1,2) } = s_{ (0,3) } = 0,\ $$

$$ s_{ (4,0) } = 1,\ s_{ (3,1) } = s_{ (2,2) } = s_{ (1,3) } = 0,\ s_{ (0,4) } =  2. $$
We shall use the following indexation for $K$:
\begin{equation}
\label{pf3_72}
\mathbf{k}_0 = (0,0),\ \mathbf{k}_1 = (1,0),\ \mathbf{k}_2 = (0,1),\ \mathbf{k}_3 = (2,0),\  \mathbf{k}_4 = (1,1),\ \mathbf{k}_5 = (0,2).    
\end{equation}
Consider the associated Hilbert space $H$ (see Introduction), and denote
\begin{equation}
\label{pf3_74}
x_j = \left[
\mathbf{t}^{ \mathbf{k}_j }
\right]_{ \mathfrak{L} },\qquad           j\in\mathbb{Z}_{0,5}.
\end{equation}
The matrix $\Gamma$ from~(\ref{f2_23}) now has the following form:
$$ \Gamma = 
\left( s_{\mathbf{k}_j + \mathbf{k}_m} \right)_{m,j=0}^5
=
\left(
(x_j,x_m)_H
\right)_{m,j=0}^5 = $$ 
\begin{equation}
\label{pf3_76}
= 
\left(
\begin{array}{cccccc} 9 & -1 & 0 & 1 & 0 & 2 \\
-1 & 1 & 0 & -1 & 0 & 0 \\ 
0 & 0 & 2 & 0 & 0 & 0 \\
1 & -1 & 0 & 1 & 0 & 0 \\
0 & 0 & 0 & 0 & 0 & 0 \\
2 & 0 & 0 & 0 & 0 & 2 \end{array}
\right).
\end{equation}
We have $\Omega_l = \Omega_0 = \{ 0,1,2 \}$, and 
\begin{equation}
\label{pf3_78}
\Gamma_l = 
\left(
\begin{array}{ccc} 9 & -1 & 0  \\
-1 & 1 & 0 \\ 
0 & 0 & 2  
\end{array}
\right),\quad l=1,2,
\end{equation}
\begin{equation}
\label{pf3_80}
\widehat\Gamma_1 = 
\left(
\begin{array}{ccc} 1 & -1 & 0  \\
-1 & 1 & 0 \\ 
0 & 0 & 0  
\end{array}
\right),\quad
\widehat\Gamma_2 = 
\left(
\begin{array}{ccc} 2 & 0 & 0  \\
0 & 0 & 0 \\ 
0 & 0 & 2  
\end{array}
\right).
\end{equation}
Conditions $\Gamma\geq 0$, and  
$$ \mathop{\rm Ker}\nolimits \Gamma_1 \subseteq \mathop{\rm Ker}\nolimits \widehat \Gamma_1,\qquad 
\mathop{\rm Ker}\nolimits \Gamma_2 \subseteq \mathop{\rm Ker}\nolimits \widehat \Gamma_2,
$$
can be verified directly.
Apply the Gram-Schmidt orthogonalization procedure to $x_0,x_1,...,x_5$.
Using~(\ref{pf3_76}) we get an orthonormal basis $\mathcal{G} = \{ g_j \}_{j=0}^3$ in $H$:
$$ g_0 = \frac{1}{3} x_0,\quad g_1 = \frac{3}{ 2\sqrt{2} } \left(
x_1 + \frac{1}{9} x_0 \right),\quad 
g_2 = \frac{1}{ \sqrt{2} } x_2,\ $$
$$ g_3 = \sqrt{ \frac{2}{3} } \left(
x_5 - \frac{1}{4} x_0 - \frac{1}{4} x_1 \right). $$
Moreover, $\{ g_j \}_{j=0}^2$ is an orthonormal basis in $H_0 = D(M_1) = D(M_2)$.
Consequently, we have $i_s = 1$.

Let $R_1$ and $R_2$ be arbitrary self-adjoint extensions of $M_1$ and $M_2$, respectively. 
We denote by
$\mathcal{R}_1$ and $\mathcal{R}_2$ the matrices of $R_1$ and $R_2$, respectively.
Using relations~(\ref{pf3_7}),(\ref{pf3_9}),(\ref{pf3_12}) and the definition of the associated operators $M_1,M_2$ we may
calculate the corresponding matrices $A_1,A_2,B_1,B_2$:
\begin{equation}
\label{pf3_82}
A_1 = 
\left(
\begin{array}{ccc} -\frac{1}{9} & \frac{ 2\sqrt{2} }{9} & 0  \\
\frac{ 2\sqrt{2} }{9} & -\frac{8}{9} & 0 \\ 
0 & 0 & 0  
\end{array}
\right),\quad
A_2 = 
\left(
\begin{array}{ccc} 0 & 0 & \frac{ \sqrt{2} }{3}  \\
0 & 0 & \frac{1}{6} \\ 
\frac{ \sqrt{2} }{3} & \frac{1}{6} & 0  
\end{array}
\right),
\end{equation}
\begin{equation}
\label{pf3_84}
B_1 = 
\left(
\begin{array}{ccc} 0 & 0 & 0  
\end{array}
\right),\quad
B_2 = 
\left(
\begin{array}{ccc} 0 & 0 & \frac{ \sqrt{3} }{2}  
\end{array}
\right).
\end{equation}
Conditions~(\ref{pf3_16}),(\ref{pf3_22}) are satisfied. Since
$$ A_1 B_2^* = A_2 B_1^* =
\left(
\begin{array}{ccc} 0 \\
0 \\ 
0\end{array}
\right),\quad B_2 B_1^* = B_1 B_2^* = 0, $$
then conditions~(\ref{pf3_18}),(\ref{pf3_20}) are satisfied with $C_1 = C_2 = 0$.
By Proposition~\ref{p3_1} we conclude that the moment problem is solvable.

\noindent
\textbf{Example 2} ($i_s=2$).
Let $n=2$, $K = K_2$, and $\mathcal{K} = K+K = K_4$.
Consider the moment problem~(\ref{f1_1}) with the following moments $\mathcal{S} = (s_{\mathbf{k}})_{\mathbf{k}\in\mathcal{K}}$:

$$ s_{ (0,0) } = 8,\ s_{ (1,0) } = s_{ (0,1) } = 0,\  s_{ (2,0) } = 2,\ s_{ (1,1) } = 0,\ s_{ (0,2) } = 2,\ $$

$$ s_{ (3,0) } =  s_{ (2,1) } = s_{ (1,2) } = s_{ (0,3) } = 0,\ $$

$$ s_{ (4,0) } = 2,\ s_{ (3,1) } = s_{ (2,2) } = s_{ (1,3) } = 0,\ s_{ (0,4) } =  2. $$
We shall use the indexation~(\ref{pf3_72}) for $K$.
Consider the associated Hilbert space $H$, and denote $x_j$ as in~(\ref{pf3_74}).
The matrix $\Gamma$ from~(\ref{f2_23}) has the following form:
$$ \Gamma = 
\left( s_{\mathbf{k}_j + \mathbf{k}_m} \right)_{m,j=0}^5
=
\left(
(x_j,x_m)_H
\right)_{m,j=0}^5 = $$ 
\begin{equation}
\label{pf3_86}
= 
\left(
\begin{array}{cccccc} 8 & 0 & 0 & 2 & 0 & 2 \\
0 & 2 & 0 & 0 & 0 & 0 \\ 
0 & 0 & 2 & 0 & 0 & 0 \\
2 & 0 & 0 & 2 & 0 & 0 \\
0 & 0 & 0 & 0 & 0 & 0 \\
2 & 0 & 0 & 0 & 0 & 2 \end{array}
\right).
\end{equation}
We have $\Omega_l = \Omega_0 = \{ 0,1,2 \}$, and 
\begin{equation}
\label{pf3_88}
\Gamma_l = 
\left(
\begin{array}{ccc} 8 & 0 & 0  \\
0 & 2 & 0 \\ 
0 & 0 & 2  
\end{array}
\right),\quad l=1,2,
\end{equation}
\begin{equation}
\label{pf3_90}
\widehat\Gamma_1 = 
\left(
\begin{array}{ccc} 2 & 0 & 0  \\
0 & 2 & 0 \\ 
0 & 0 & 0  
\end{array}
\right),\quad
\widehat\Gamma_2 = 
\left(
\begin{array}{ccc} 2 & 0 & 0  \\
0 & 0 & 0 \\ 
0 & 0 & 2  
\end{array}
\right).
\end{equation}
Conditions $\Gamma\geq 0$, and  
$$ \mathop{\rm Ker}\nolimits \Gamma_1 \subseteq \mathop{\rm Ker}\nolimits \widehat \Gamma_1,\qquad 
\mathop{\rm Ker}\nolimits \Gamma_2 \subseteq \mathop{\rm Ker}\nolimits \widehat \Gamma_2,
$$
can be checked directly.
Apply the Gram-Schmidt orthogonalization procedure to $x_0,x_1,...,x_5$.
We get an orthonormal basis $\mathcal{G} = \{ g_j \}_{j=0}^4$ in $H$:
$$ g_0 = \frac{1}{ 2\sqrt{2} } x_0,\quad g_1 = \frac{1}{ \sqrt{2} } x_1,\quad 
g_2 = \frac{1}{ \sqrt{2} } x_2,\quad 
g_3 = \sqrt{ \frac{2}{3} } \left(
x_3 - \frac{1}{4} x_0 \right), $$
$$ g_4 = \frac{ \sqrt{3} }{2} \left(
x_5 + \frac{1}{3} x_3 - \frac{1}{3} x_0 \right). $$
Observe that $\{ g_j \}_{j=0}^2$ is an orthonormal basis in $H_0 = D(M_1) = D(M_2)$.
Therefore $i_s = 2$.

Let $R_1$ and $R_2$ be arbitrary self-adjoint extensions of $M_1$ and $M_2$, respectively. 
Denote by
$\mathcal{R}_1$ and $\mathcal{R}_2$ the matrices of $R_1$ and $R_2$, respectively.
By relations~(\ref{pf3_7}),(\ref{pf3_9}),(\ref{pf3_12}) we may
calculate the corresponding matrices $A_1,A_2,B_1,B_2$:
\begin{equation}
\label{pf3_92}
A_1 = 
\left(
\begin{array}{ccc} 0 & \frac{ 1 }{2} & 0  \\
\frac{ 1 }{2} & 0 & 0 \\ 
0 & 0 & 0  
\end{array}
\right),\quad
A_2 = 
\left(
\begin{array}{ccc} 0 & 0 & \frac{1}{2}  \\
0 & 0 & 0 \\ 
\frac{1}{2} & 0 & 0  
\end{array}
\right),
\end{equation}
\begin{equation}
\label{pf3_94}
B_1 = 
\left(
\begin{array}{ccc} 0 & \frac{ \sqrt{3} }{2} & 0 \\
0 & 0 & 0  
\end{array}
\right),\quad
B_2 = 
\left(
\begin{array}{ccc} 0 & 0 & \frac{ -1 }{ 2 \sqrt{3} } \\
0 & 0 & \frac{1}{2}
\end{array}
\right).
\end{equation}
Conditions~(\ref{pf3_16}),(\ref{pf3_22}) are fulfilled.
Since
$$ A_1 B_2^* = A_2 B_1^* =
\left(
\begin{array}{ccc} 0 & 0 \\
0 & 0 \\ 
0 & 0 \end{array}
\right),\quad B_2 B_1^* = B_1 B_2^* = 
\left(
\begin{array}{cc} 0 & 0 \\
0 & 0 \end{array}
\right), $$
then conditions~(\ref{pf3_18}),(\ref{pf3_20}) are satisfied with 
$$ C_1 = C_2 = \left(
\begin{array}{cc} 0 & 0 \\
0 & 0 \end{array}
\right). $$
By Proposition~\ref{p3_1} we obtain that this moment problem is solvable.

}

\noindent
Address:

V. N. Karazin Kharkiv National University \newline\indent
School of Mathematics and Computer Sciences \newline\indent
Department of Higher Mathematics and Informatics \newline\indent
Svobody Square 4, 61022, Kharkiv, Ukraine

Sergey.M.Zagorodnyuk@gmail.com; zagorodnyuk@karazin.ua

\end{document}